\documentclass[journal, 10pt, conference]{ieeeconf}      

\IEEEoverridecommandlockouts                              
\usepackage{graphics} 
\usepackage{epsfig} 
\usepackage{mathptmx} 
\usepackage{times} 
\usepackage{amsmath} 
\usepackage{amssymb}  
\usepackage{subfigure}
\usepackage{wasysym}

\usepackage{color}
\newtheorem{definition}{Definition}
\newtheorem{theorem}{Theorem}
\newtheorem{lemma}{Lemma}

\title{\LARGE \bf Optimal Regional Tracking Control of Time-Fractional Diffusion Systems}

\author{Fudong Ge$^{1}$ and YangQuan Chen$^{2}$
\thanks{This work was accepted in the 2021 Online American Control Conference (ACC2021).}
\thanks{$^{1}$School of Computer Science, China University of Geosciences, Wuhan 430074, PR China
{\tt\small (Email: gefd@cug.edu.cn)}}
\thanks{$^{2}$School of Engineering (MESA-Lab), University of California, Merced, CA 95343, USA
{\tt\small (Email: ychen53@ucmerced.edu)}}
}

\begin{document}

\maketitle
\thispagestyle{empty}
\pagestyle{empty}

\begin{abstract}
In this paper, we aim to explore optimal regional trajectory tracking control problems of the anomalous subdiffusion processes
governed by time-fractional diffusion systems under the Neumann boundary conditions. Using eigenvalue theory of the system
operator and the semigroup theory, we explore the existence and some estimates of the mild solution to the considered system.
An approach on finding solution to the optimal problem that minimizes the regional trajectory tracking error and the corresponding
control cost over a finite space and time domain is then explored via the Hilbert uniqueness method (HUM). The obtained results
not only can be directly used to investigate the systems that are not controllable on the whole domain, but also yield an explicit
expression of the control signal in terms of the desired trajectory. Most importantly, it is worth noting that our results in this
paper are still novel even for the special case when the order of fractional derivative is equal to one. Finally, we provide a
numerical example to illustrate our theoretical results.
\end{abstract}
\label{sec1}

\section{Introduction}
\noindent

During the past two decades, there is an increasing activity in the discussion of tracking problems for conventional reaction-diffusion dynamic systems, which can be
divided into two steps: trajectory planning and tracking control (see e.g., the monographes \cite{L2017Optimal,meurer2012control}). In these studies, the trajectory
planning step attempts to generate a reference trajectory for the given desired function, while the tracking control focuses on system dynamics and hopes to
design a sequence of inputs to track the pre-planned reference trajectory \cite{Choi2004PID}. For the trajectory planning problems, we refer the reader to
\cite{meurer2009trajectory,meurer2011flatness,freudenthaler2020pde-based} where the flatness-based feedforward control strategies were presented or to
\cite{2004Diffusion} where the hardware and numerical illustrations were carried out. Moreover, to improve the accuracy in tracking control, various controller design techniques such as sliding-mode control \cite{pisano2011tracking},
robust control \cite{weng2017data-driven}, iterative learning control \cite{Kang2014A} and active disturbance rejection control (ADRC)
\cite{Zhou2019TAC} have been developed.

On the other hand, after the pioneering work given by Einstein in \cite{einstein2005uber}, it is confirmed that conventional diffusion system can well model
the Brownian motion, whose mean-square displacement (MSD) is a linear function of time $t$. However, there exist a great deal of extremely complex transport
processes that
is characterized by a power-law MSD relation $($i.e., MSD$=t^\alpha, \alpha>0)$ including the anomalous subdiffusion case with $\alpha\in (0,1)$ and the
anomalous superdiffusion case
with $\alpha>1$. In these anomalous situations, the usual physical laws would never be followed and the mathematical models will divert from the traditional
integer-order systems to the fractional-order cases \cite{Ge2019InforSci,Ge2019ND,uchaikin2012fractional,Ge2018Book}.
Besides,  we see that for anomalous subdiffusion process, time-fractional diffusion system has been confirmed as a powerful tool to model it
\cite{ralf,Ge2016FCAA,Ge2018SCL}.
Here the time-fractional diffusion system is a new extension model of conventional diffusion system by replacing the first order time derivative
with a fractional-order derivative of order $\alpha\in(0,1]$. This is due to the fact that fractional-order derivative is defined as a kind of
convolution hence representing well the dynamics inheriting subdiffusive properties and moreover, the fractional-order derivative would recover
the first-order derivative if it approaches to one \cite{kilbas}.
Then, based on our previous work on trajectory planning problem of time-fractional reaction-diffusion systems \cite{Ge2020IFAC}, in this paper,
we go on investigating the optimal trajectory tracking control problems of linear time-fractional diffusion systems with the Neumann boundary
conditions. Further results on optimal tracking control of coupled nonlinear time-fractional diffusion systems under more general boundary
conditions will be discussed in our forthcoming works.

Let $\Omega$ be an open bounded subset of $\mathbf{R}^n$ with Lipschitz continuous boundary $\partial\Omega$
and denote $Q=\Omega\times \left[0,T\right]$, $\Sigma= \partial \Omega \times[0,T]$ with $T>0.$ Herein, we consider the following time-fractional diffusion
system with a Caputo fractional derivative ${}_0^CD^{\alpha}_t$ of order $\alpha\in (0,1]$:
\begin{equation}\label{problem}
\left\{{\begin{array}{*{20}{l}}
{{}_0^CD^{\alpha}_{t}y(x,t)=Ay(x,t)+u(x,t)\mbox{ in }Q,}\\
{\frac{\partial y}{\partial \nu}(x,t)=0 \mbox{ in } \Sigma,}\\
{ y(x,0)=y_0(x)\mbox{ in } \Omega,}
\end{array}}\right.
\end{equation}
where $A$ generates a strongly continuous semigroup $\{\Phi(t)\}_{t\geq 0}$ on the Hilbert space $L^2(\Omega)$, $-A$ is a uniformly elliptic operator
(see e.g., the Definition 9.2 of \cite{Renardy2016an}), $u\in L^2(Q)$ denotes the control inputs and $\nu$ represents the unit outside normal vector of the boundary
$\partial \Omega$. Here $L^2(\Omega)$ represents the usual Hilbert space endowed with the inner product $(\cdot,\cdot)_{L^2(\Omega)}$ and the norm
$\|\cdot\|_{L^2(\Omega)}$.
As cited in \cite{sun2018new}, system $(\ref{problem})$ covers a great deal of real-world applications in a spatially inhomogeneous environment.
Typical examples include the reheating processes of heterogeneous metal slabs \cite{Ge2018Book} or the flow through porous media with varying sources
or sinks \cite{uchaikin2008fractional} and so on.



Taking into account that not all the states of time-fractional diffusion systems are reachable in the whole domain of interest. To address
this issue, regional control ideas such as regional controllability \cite{Ge2016Automatica2,Ge2016JMAA}, regional observability \cite{Ge2016Automatica1}
and regional stability \cite{Kang2020Constrained,Ge2020IJRNC} have been employed and well studied due to their advantages of offering potential to
reduce computational requirements and being possible to study the systems that are not controllable on the whole domain. With these in mind, the
novelty in this paper is promoting to study regional trajectory tracking control problem of the system $(\ref{problem})$. For this purpose,
we focus on employing optimal control strategy to determine control signals by minimizing the proposed tracking cost functional. However, optimal
control design for system $(\ref{problem})$ is very challenging or even impossible due to the infinite dimensionality property of the problem.
To overcome this limitation, the HUM provides an alternative approach \cite{Ge2016Automatica2,Lions1971,glowinski2008exact}. In this method, dual
system is selected to determine the explicit expression of optimal solution to the cost functional. To the best of our knowledge,
no results are available on this topic. Most importantly, we claim that our results in this
paper are still novel even for the special case when the order of fractional derivative in considered system is equal to one.

The rest of this paper is organized as follows. Some preliminary results that are useful for the study are
given in Section 2. In Section 3, we present our main results on the HUM-based optimal regional trajectory  tracking control
strategy for time-fractional diffusion systems. This is illustrated in Section 4, where
a numerical example is presented.

\section{Preliminaries}

\begin{definition}\cite{kilbas} Given a function $\phi: [0,T]\to \mathbf{R}$,
the Riemann-Liouville fractional integral of order $\alpha>0$ for $\phi$ is as follows
\begin{eqnarray}\label{fractionainte}
_0I^{\alpha}_t \phi(t)= \int^{t}_{0}{\frac{(t-s)^{\alpha-1}}{\Gamma(\alpha) }\phi(s)ds}.
\end{eqnarray}
Here
\begin{eqnarray}
\Gamma(\alpha)=\int_0^\infty{t^{\alpha-1}e^{-t}}dt\end{eqnarray}
denotes the Euler gamma function and the right side of $(\ref{fractionainte})$ is pointwise defined on $[0,T]$.
\end{definition}
\begin{definition}\cite{kilbas}\label{fderivative}
Given a function $\phi: [0,T]\to \mathbf{R}$,
the Caputo fractional derivative of order $\alpha\in (0,1]$ for $\phi$ is
\begin{eqnarray}
^C_0D_t^\alpha \phi(t)=\left\{\begin{array}{l}{}_0I_t^{1-\alpha} \frac{\partial}{\partial t}\phi(t),~~0< \alpha <1,
\\
\frac{\partial}{\partial t}\phi(t),~~~~~~~~~~\alpha =1
\end{array}
\right.
\end{eqnarray}
provided that the right side is pointwise defined on $[0,T]$.
\end{definition}

Recall that $-A$ is a uniformly elliptic operator, by \cite{Hilbert}, under the Neumann boundary condition, the eigenvalue
pairing $(\lambda_k,\xi_k)_{k\in \mathbf{N}^n}$ with $\mathbf{N}=\{0,1,2,\cdots\}$ of the Sturm-Liouville problem
\begin{eqnarray}
\begin{array}{l}
A\xi_k=\lambda_k \xi_k, ~k=(k_1,k_2,\cdots,k_n)\in \mathbf{N}^n
\end{array}
\end{eqnarray}
satisfies
\begin{eqnarray*}
\begin{array}{l}
\lambda_0=0,~\lambda_k<0 \mbox{ for all }k\in \mathbf{N}^n\backslash \{0\} \mbox{ with } \lim\limits_{k_i\to \infty }\lambda_k=-\infty
\end{array}
\end{eqnarray*}
and $\{\xi_k(x) \}_{k\in \mathbf{N}^n}$ forms a orthonormal basis of $ L^2(\Omega)$. This is general. For example, when
$\Omega=\{x=(x_1,\cdots,x_n)\in \mathbf{R}^n: 0<x_i<L_i, ~i\in I_n=\{1,2,\cdots,n\} \}$ is a $1\leqslant n-$dimensional parallelepipedon, we
refer the reader to Lemma 2 of \cite{meurer2011flatness} for a detailed expression of the corresponding eigenvalue pairing $(\lambda_k,\xi_k)_{k\in \mathbf{N}^n}$.
With these,  any $\varphi\in L^2(\Omega)$ can be expressed as
\begin{eqnarray}
\begin{array}{l}
\varphi(x)=\sum\limits_{k\in \mathbf{N}^n}{\left(\varphi,\xi_k\right)_{L^2(\Omega)} \xi_k(x)}.
\end{array}
\end{eqnarray}
Then, the strongly continuous semigroup $\{\Phi(t)\}_{t\geq 0}$ on  $L^2(\Omega)$ generated by $A$ satisfies
\begin{eqnarray}
\begin{array}{l}
\Phi(t) \varphi=(\varphi,\xi_0)_{L^2(\Omega)}\xi_0+\sum\limits_{k\in \mathbf{N}^n\backslash \{0\}}{e^{\lambda_kt}
\left(\varphi,\xi_k \right)_{L^2(\Omega)}\xi_k}\\{\kern 26 pt}
=\sum\limits_{k\in \mathbf{N}^n}{e^{\lambda_kt}\left(\varphi,\xi_k\right)_{L^2(\Omega)} \xi_k}.
\end{array}
\end{eqnarray}

Denote by $H^2(\Omega) $ and $ H_0^1(\Omega)$ the usual Sobolev spaces (see e.g., \cite{adams2003sobolev}), now we are ready to give
the following lemma.

\begin{lemma}\cite{Ge2020IJRNC,japan}\label{solution}
If $y_0\in L^2(\Omega)$, then there exists a unique mild
solution $y\in C\left([0,T];L^2(\Omega)\right)\cap C\left((0,T];H^2(\Omega)\cap H_0^1(\Omega)\right)$ to system $(\ref{problem})$ such that
${}_0^CD^\alpha_t y\in C\left((0,T];L^2(\Omega)\right)$ and the estimate
\begin{eqnarray}
\begin{array}{l}
\|y\|_{C\left([0,T];L^2(\Omega)\right)}
\leqslant \gamma_1 \|y_0\|_{L^2(\Omega)} +\gamma_2\left\|u\right\|_{L^2\left(Q\right)}
\end{array}
\end{eqnarray}
holds true for some $\gamma_1,\gamma_2>0$. Moreover, $y(x,t)$ satisfies
\begin{eqnarray}
\begin{array}{l}
y(x,t)=\mathcal{M}_{\alpha}(t)y_0(x)
+\int_0^t{ \frac{\mathcal{K}_{\alpha}(t-\tau)}{(t-\tau)^{1-\alpha}} u(x,\tau)}d\tau,
\end{array}
\end{eqnarray}
where
\begin{eqnarray}
\begin{array}{l}
\mathcal{M}_{\alpha}(t)\varphi(x)= \sum\limits_{k\in \mathbf{N}^n} E_{\alpha}(\lambda_k t^\alpha)(\varphi,\xi_k)_{L^2(\Omega)}\xi_k(x),
\end{array}
\end{eqnarray}
\begin{eqnarray}
\begin{array}{l}
\mathcal{K}_{\alpha}(t)\varphi(x)=\sum\limits_{k\in \mathbf{N}^n}E_{\alpha,\alpha} (\lambda_kt^\alpha)\left(\varphi,\xi_k\right)_{L^2(\Omega)}
\xi_k(x)
\end{array}
\end{eqnarray}
with $ \varphi\in L^2(\Omega)$ and
\begin{eqnarray}\label{eq2.2}
\begin{array}{l}
E_{\alpha,\beta}(t)=\sum\limits_{k=0}^\infty{\frac{t^k}{\Gamma(\alpha k+\beta)}},~\alpha, \beta>0, ~t\geqslant  0,
\end{array}
\end{eqnarray}
denotes the Mittag-Leffler function in two parameters. In particular, we write $E_{\alpha}(t)=E_{\alpha,1}(t)$ for short when $\beta=1$.
\end{lemma}

Let $y(x,t,u)$ denote the solution of system $(\ref{problem})$ for any given control $u\in L^2(Q)$.
Choose $\omega\subseteq \Omega$  a positive Lebesgue measure sub-region, we define $\left(L^2(\omega),\|\cdot\|_{L^2(\omega)}\right)$
as the corresponding Hilbert space on $\omega$ and set $Q_\omega=\omega\times [0,T]$. Moreover, let $\chi_\omega$ be the characteristic
function of $\omega$ given by
\begin{eqnarray}
\chi_\omega(x)=\left\{\begin{array}{l}
1,~\mbox{ if } x\in \omega,\\
0,~\mbox{ if } x\in \Omega\backslash \omega.
\end{array}\right.
\end{eqnarray}
We have the following definition.
\begin{definition}
The considered optimal regional trajectory tracking control problem for system $(\ref{problem})$ in $\omega$ at time $T$ concerns
how to design controller $u\in L^2(Q)$ such that $y(x,t,u)$ starting from $y_0(x)$ could reach the given target
function $y_{dT}=y_d(\cdot,T)\in L^2(\omega)$ in $\omega$ as close as possible along with the given trajectory $y_d(x,t)\in L^2\left(Q_\omega\right)$
within $t\in [0,T]$, i.e.,
\begin{eqnarray}\label{probformu}
\begin{array}{l}
~~~~~~y_0(x)\xrightarrow [{\mbox{along with } y_d(x,t)}]{\mbox{optimal } u \in L^2(Q)\mbox{ in } y(x,t,u)}y_{dT}(x), ~x\in \overline{\omega}.
\end{array}
\end{eqnarray}
\end{definition}
\vspace{0.2cm}

The following lemma plays a key role for our later technical development.

\begin{lemma}\label{integralbypart}\cite{Ge2021ISATrans}
Let
\begin{eqnarray}\label{RLintegral}
{}_tI_T^\alpha \varphi(t)=\frac {1}{\Gamma(\alpha)} \int_t^T(s-t)^{\alpha-1}\varphi(s)ds,~\alpha \in (0,1]
\end{eqnarray}
be the right-sided Riemann-Liouville fractional integral and
let us denote ${}_tD^{\alpha}_{T}$ the right-sided Riemann-Liouville fractional derivative given by (see Section 2.1 of \cite{kilbas})
\begin{eqnarray}\label{RLCaputo}
{}_tD_T^\alpha \varphi(t)=
\left\{\begin{array}{l}
- \frac{d}{dt} {}_tI_T^{1-\alpha} \varphi(t),~\alpha\in (0,1),\\
-\frac{d }{d t}\varphi(t),~\alpha=1.
\end{array}
\right.
\end{eqnarray}
For any $\alpha\in (0,1]$, if $\varphi_1'\in L^p(0,T)$, $\varphi_2\in L^q(0,T)$, $p,q\geqslant 1$, $1/p+1/q\leqslant 1+\alpha$ and
$p\neq 1$, $q\neq1$ in the case when  $1/p+1/q= 1+\alpha$, then the formula
\begin{eqnarray}\label{formubypart}
\begin{array}{l}
\int_0^T{\varphi_2(t){}^C_0D_t^{\alpha}\varphi_1(t)}dt
=\int_0^T{\varphi_1(t){}_tD_T^{\alpha}\varphi_2(t)}dt\\{\kern 95 pt}
+\left[\varphi_1(t){}_tI_T^{1-\alpha}\varphi_2(t)\right]_{t=0}^{t=T}
\end{array}
\end{eqnarray}
holds true.
\end{lemma}

\section{Optimal regional tracking control}

The control objective of this article is the guidance of state trajectories of system $(\ref{problem})$ along a
given trajectory $y_d\in L^2(Q_\omega)$ to reach the final target function $y_{dT}=y_d(\cdot,T)\in L^2(\omega)$
in $\omega$ at time $T$. The closer the controlled $y(x,t)$, $y(x,T)$ follows the desired trajectory $y_d(x,t)$
and the target $y_{dT}(x)$, the better the control target is achieved.

Suppose that $U_{ad} $ is an nonempty closed, convex subset of $L^2(Q)$, in this section, we aim to find $u_r\in U_{ad}$
such that
\begin{eqnarray}\label{minimization}
\begin{array}{l}
J(u_r)\leqslant J(u)\mbox{ for all } u\in U_{ad},
\end{array}
\end{eqnarray}
where $J$ is the squared difference integrated performance functional given by
\begin{eqnarray}\label{trackingcostf}
\begin{array}{l}
J(u)=\frac{r_1}{2} \int_{Q_\omega} {\left|\chi_\omega y(x,t,u)- y_d(x,t)\right|^2}dxdt\\{\kern 29 pt}
+\frac{r_2}{2}\int_{\omega} {\left|\chi_\omega y(x,T,u)-y_{dT}(x)\right|^2}dx
\\{\kern 29 pt}
+ \frac{r_3}{2}\int_{Q}{ \left|u(x,t)\right|^2} dxdt
\end{array}
\end{eqnarray}
and $r_1,$$r_2\geqslant 0$, $r_3>0$ are three given constants. Here $y_{dT}$ may be the equilibrium or any
given target functions of system $(\ref{problem})$.
For this purpose, since $U_{ad}\subseteq L^2(Q)$ is a nonempty closed convex subset, the following lemma is necessary.
\begin{lemma}\cite{Lions1971,glowinski2008exact}\label{Lions1971}
 Assume that the quadratic function $u\to J(u)$ is strictly convex and differentiable that  satisfies
\begin{equation}
J(u)\to \infty~\mbox{ as } \|u\|\to \infty, ~u\in U_b.
\end{equation}
Then the uniqueness element $u$ in closed, convex subset $U_{ad}$ satisfying $J(u)=\inf\limits_{v\in U_{ad}}{J(v)}$ is characterized by
\begin{equation}
J'(u)\cdot(v-u)\geqslant  0,~~\forall v\in U_{ad}.
\end{equation}
\end{lemma}

Based on Lemma $\ref{Lions1971}$,
we get that the unique solution of the optimization problem $(\ref{minimization})$
can be characterized by
\begin{eqnarray}\label{eq34}
\begin{array}{l}
J'(u_r)\cdot(u-u_r)\geqslant 0\mbox{ holds true for all }u\in U_{ad}
\end{array}
\end{eqnarray}
and more precisely,
\begin{equation}\label{OptimalCond}
\begin{array}{l}
r_1\int_{Q} {\left( p_\omega y(x,t,u_r)- \chi_\omega^* y_d(x,t)\right)\left( y(x,t,u)-  y(x,t,u_r)\right)}dxdt\\
+r_2\int_{\Omega}{\left( p_\omega y(x,T,u_r)-\chi_\omega^* y_{dT}(x) \right)\left( y(x,T,u)- y(x,T,u_r) \right)}dx\\
+ r_3\int_Q{u_r(x,t)\left(u(x,t)-u_r(x,t)\right)}dx dt \geqslant 0~\mbox{ for all }u\in U_{ad}
\end{array}
\end{equation}
 after a simple duality derivation. Here
\begin{eqnarray}
\begin{array}{l}
\chi_\omega^* \varphi(x)=\left\{\begin{array}{l}\varphi(x), ~~x\in \omega ,
 \\0,~~~x\in \Omega\backslash \omega
\end{array}\right.
\end{array}
\end{eqnarray}
denotes the adjoint operator of $\chi_\omega$ and $p_\omega=\chi_\omega^* \chi_\omega $.
Further, to simplify above equation $(\ref{OptimalCond})$, let us introduce the following adjoint system
\begin{equation}\label{observer}
\left\{{\begin{array}{*{20}{l}}
{}_tD^{\alpha}_{T}z(x,t)=A^* z(x,t)
+r_1\left(p_\omega y(x,t,u_r)-\chi_\omega^*  y_d(x,t)\right)\\{\kern 48 pt}
\mbox{ in }Q,\\
\frac{\partial z}{\partial \nu}(x,t)=0
 \mbox{ in } \Sigma,\\
\lim\limits_{t\to 0^+}{}_0I_t^{1-\alpha} Rz(x,t)=r_2(p_\omega  y(x,T,u_r)-\chi_\omega^* y_{dT}(x))\\{\kern 49 pt}
\mbox{ in } [0,T],
\end{array}}\right.
\end{equation}
where $A^*$ is the adjoint operator of $A$ and $R$ is an operator given by $Rz_1(t)=z_1(T-t)$ as in property 2.7 of \cite{klimek2009solutions}
satisfying
\begin{equation}\label{Rintegral}
\begin{array}{l}
R\left({}_tI_T^\alpha z_1(t)\right)=\frac {1}{\Gamma(\alpha)} \int_{T-t}^T(s-T+t)^{\alpha-1}z_1(s)ds\\{\kern 47 pt}
=\frac {1}{\Gamma(\alpha)} \int_0^t(t-\tau)^{\alpha-1}Rz_1(\tau)d\tau\\{\kern 47 pt}
={}_0I_t^\alpha Rz_1(t),
\end{array}
\end{equation}
\begin{equation}\label{Rderivative}
\begin{array}{l}
R\left({}_tD_T^\alpha z_1(t)\right)
=\frac {-\frac{d}
{d (T-t)}}{\Gamma(1-\alpha)} \int_{T-t}^T(s-T+t)^{-\alpha}z_1(s)ds\\{\kern 50 pt}
=\frac {1}{\Gamma(1-\alpha)}\frac{d}{d t}\int_0^t(t-\tau)^{-\alpha}z_1(T-\tau)d\tau\\{\kern 50 pt}
={}_0D_t^\alpha Rz_1(t)
\end{array}
\end{equation}
as a consequence of Eq.$(\ref{RLintegral})$ and Eq.$(\ref{RLCaputo})$.
To establish the existence of a unique solution to system $(\ref{observer})$, it is supposed that the eigenvalue pairing of operator $A^*$ under
the Neumann boundary conditions is $\left(\lambda_k^*,\xi_k^*\right)_{k\in \mathbf{N}^n}$, where $\{\xi_k^*(x) \}_{k\in \mathbf{N}^n}$ also forms a orthonormal
basis of $ L^2(\Omega)$.
Recall from Lemma $\ref{solution}$ and the Lemma 1 of \cite{Ge2016Automatica2}, given any $y_d\in L^2\left(Q\right)$ and $y_{dT}\in L^2(\Omega)$,
based on the property of
operator $R$, we have the following result.

\begin{lemma}\label{backsolution}
Given any $y_d\in L^2\left(Q\right)$ and $y_{dT}\in L^2(\Omega)$, if conditions of Lemma $\ref{solution}$ are satisfied, then
system $(\ref{observer})$ admits a unique mild solution $z\in C\left([0,T];L^2(\Omega)\right)\cap C\left((0,T];H^2(\Omega)\cap H_0^1(\Omega)\right)$ as follows
\begin{eqnarray}\label{seriesexpansion}
\begin{array}{l}
z(x,t,u_r)=r_2\frac{\mathcal{K}_{\alpha}^*(T-t)}{(T-t)^{1-\alpha}}
\left(p_\omega y(x,T,u_r)-\chi_\omega^*  y_{dT}(x)\right)\\{\kern 10pt}
+r_1\int_t^T{ \frac{\mathcal{K}_{\alpha}^*(\varsigma-t)}{(\varsigma-t)^{1-\alpha}} \left(p_\omega y(x,\varsigma,u_r)-\chi_\omega^*
y_d(x,\varsigma)\right)}d\varsigma,
\end{array}
\end{eqnarray}
where
\begin{eqnarray}
\begin{array}{l}
\mathcal{K}_{\alpha}^*(t)\varphi(x)=\sum\limits_{k\in \mathbf{N}^n}E_{\alpha,\alpha} (\lambda_k^* t^\alpha)
\left(\varphi,\xi_k^*\right)_{L^2(\Omega)}\xi_k^*(x)
\end{array}
\end{eqnarray}
for any $\varphi\in L^2(\Omega)$.
\end{lemma}
\begin{proof} Taking $R-$transformation defined above on both sides of system $(\ref{observer})$, using Eq.$(\ref{Rderivative})$, one has
\begin{eqnarray*}\label{observertransfer}
\left\{{\begin{array}{*{20}{l}}
{}_0D_t^\alpha R z(x,t)=A^*Rz(x,t)\\{\kern 41 pt}
+r_1\left( p_\omega y(x,T-t,u_r)-\chi_\omega^*  y_d(x,T-t)\right)\mbox{ in }Q,\\
\frac{\partial}{\partial \nu}R z(x,t)=0~
 \mbox{in } \Sigma,\\
\lim\limits_{t\to 0^+}{}_0I_t^{1-\alpha} Rz(x,t)
=r_2(p_\omega y(x,T,u_r)-\chi_\omega^*  y_{dT}(x)) \\{\kern 41 pt}
\mbox{ in } [0,T].
\end{array}}\right.
\end{eqnarray*}
It follows that \cite{Ge2016Automatica2} and Lemma $\ref{solution}$ that a unique mild solution $z\in C\left([0,T];L^2(\Omega)\right)
\cap C\left((0,T];H^2(\Omega)\cap H_0^1(\Omega)\right)$ to system $(\ref{observer})$ exists and moreover, it satisfies
\begin{eqnarray*}
\begin{array}{l}
z(x,t)=R r_2 \frac{\mathcal{K}_{\alpha}^*(t)}{t^{1-\alpha}}
\left(p_\omega y(x,T,u_r)-\chi_\omega^* y_{dT}(x)\right)\\{\kern 12 pt}
+R r_1\int_0^t{ \frac{\mathcal{K}_{\alpha}^*(t-\tau)}{(t-\tau)^{1-\alpha}}  \left(p_\omega y(x,T-\tau,u_r)-\chi_\omega^*  y_d(x,T-\tau)\right)}d\tau,
\\
=r_2\frac{\mathcal{K}_{\alpha}^*(T-\tau)}{(T-t)^{1-\alpha}}
\left(p_\omega y(x,T,u_r)-\chi_\omega^*  y_{dT}(x)\right)\\{\kern 12 pt}
+r_1\int_t^T{ \frac{\mathcal{K}_{\alpha}^*(\varsigma-t)}{(\varsigma-t)^{1-\alpha}} \left(p_\omega y(x,\varsigma,u_r)-\chi_\omega^*
y_d(x,\varsigma)\right)}d\varsigma.
\end{array}
\end{eqnarray*}
This finishes the proof.
\end{proof}

In what follows, we proceed to simplify $(\ref{OptimalCond})$ based on the adjoint system  $(\ref{observer})$.

Indeed, for the first term of above equation $(\ref{OptimalCond})$, using Lemma $\ref{integralbypart}$, it yields that
\begin{eqnarray*}
\begin{array}{l}
r_1\int_Q {\left(p_\omega y(x,t,u_r)-\chi_\omega^* y_d(x,t)\right) \left(y(x,t,u)-y(x,t,u_r)\right)}dxdt\\
=\int_Q{\left({}_tD^{\alpha}_{T}z(x,t,u_r)-A^* z(x,t,u_r)\right)\left(y(x,t,u)-y(x,t,u_r)\right)}dxdt\\
=\int_Q {z^T(x,t,u_r){}^C_0D_t^{\alpha}\left(y(x,t,u)-y(x,t,u_r)\right)}dxdt\\
- \int_Q{A^* z(x,t,u_r)\left(y(x,t,u)-y(x,t,u_r)\right)}dxdt\\
-r_2\int_{\Omega} {(p_\omega y(x,T,u_r)-\chi_\omega^*  y_{dT}(x))\left(y(x,T,u)-y(x,T,u_r)\right)}dx.
\end{array}
\end{eqnarray*}
Since $-A$ is a uniformly elliptic operator, for any $\varphi_1,$ $ \varphi_2\in L^2(\Omega),$ one has \cite{Renardy2016an}
\begin{eqnarray}\label{Green}
\begin{array}{l}
\int_\Omega{\left[\varphi_1(x)A^*\varphi_2(x)-\varphi_2(x) A\varphi_1(x)\right]}dx\\
=\int_{\partial\Omega}{\left[\varphi_1(x)\frac{\partial \varphi_2(x)}{\partial \nu}-\varphi_2(x) \frac{\partial \varphi_1(x)}{\partial \nu}\right]}dx.
\end{array}
\end{eqnarray}
Then, using formula $(\ref{Green})$, we have
\begin{eqnarray*}
\begin{array}{l}
r_1\int_Q {\left(p_\omega y(x,t,u_r)-\chi_\omega^*  y_d(x,t)\right) \left(y(x,t,u)-y(x,t,u_r)\right)}dxdt\\
=\int_Q {z(x,t,u_r)  {}^C_0D_t^{\alpha}\left(y(x,t,u)-y(x,t,u_r)\right)}dxdt\\
-\int_Q {z(x,t,u_r)A\left(y(x,t,u)-y(x,t,u_r)\right)}dxdt\\
-\int_\Sigma {z(x,t,u_r)\left(\frac{\partial y(x,t,u)}{\partial \nu}-\frac{\partial y(x,t,u_r)}{\partial \nu}\right)}dxdt\\
+\int_\Sigma {\left\{\frac{\partial z (x,t,u_r)}{\partial \nu}\right\}\left(y(x,t,u)-y(x,t,u_r)\right)}dxdt\\
-r_2\int_{\Omega} {(p_\omega y(x,T,u_r)-\chi_\omega^*   y_{dT}(x))\left(y(x,T,u)-y(x,T,u_r)\right)}dx.
\end{array}
\end{eqnarray*}
This, together with the boundary conditions of considered systems, yields that
\begin{eqnarray*}
\begin{array}{l}
 r_1\int_Q {\left(p_\omega y(x,t,u_r)-\chi_\omega^*  y_d(x,t)\right) \left(y(x,t,u)-y(x,t,u_r)\right)}dxdt\\
+r_2\int_{\Omega} {(p_\omega y(\cdot,T,u)-\chi_\omega^* y_{dT}(x))\left(y(x,T,u)-y(x,T,u_r)\right)}dx\\
=\int_Q {z(x,t,u_r)  \left({}^C_0D_t^{\alpha} -A \right)\left(y(x,t,u)-y(x,t,u_r)\right)}dxdt\\
=\int_Q{z(x,t,u_r)\left(u(x,t)-u_r(x,t)\right)}dxdt.
\end{array}
\end{eqnarray*}
Therefore, the optimality condition $(\ref{OptimalCond})$ can be simplified to
\begin{eqnarray}\label{finaloptimalcond}
\begin{array}{l}
\int_Q{ \left(r_3u_r(x,t)+ z(x,t,u_r)\right) \left(u(x,t)-u_r(x,t)\right)}dx dt \geqslant 0
\end{array}
\end{eqnarray}
for all $u\in U_{ad}.$

Now we summarize the following result and omit the detailed proof.

\begin{theorem}\label{theorem4.1}
Given any target trajectory $y_d\in L^2\left(Q\right)$ and the target function $y_{dT}=y(\cdot,T)\in L^2(\Omega)$, the
optimal problem $(\ref{minimization})$ admits a unique optimal solution $u_r$ that is determined by the system $(\ref{problem})$ and
the adjoint system $(\ref{observer})$ satisfying the variational inequality $(\ref{finaloptimalcond})$.
\end{theorem}

In particular, when $U_{ad}=L^2(Q)$, since $r_3>0,$ $(\ref{finaloptimalcond})$ holds true if
\begin{eqnarray}\label{u_r}
\begin{array}{l}
u_r(x,t)= -\frac{1}{r_3}z(x,t,u_r) \mbox{ for all } (x,t)\in Q.
\end{array}
\end{eqnarray}
Furthermore, in order to obtain the optimal control $(\ref{minimization})$ in feedback form,
according to Lemma $\ref{solution}$, $(\ref{seriesexpansion})$ and $(\ref{u_r})$, we have
\begin{eqnarray*}
\begin{array}{l}
u_r(x,t)=-\frac{1}{r_3}z(x,t,u_r)\\{\kern 28 pt}
=-\frac{r_2}{r_3}\frac{\mathcal{K}_{\alpha}^*(T-t)}{(T-t)^{1-\alpha}}
\left(p_\omega y(x,T,u_r)-\chi_\omega^*  y_{dT}(x)\right)\\{\kern 39 pt}
-\frac{r_1}{r_3}\int_t^T{ \frac{\mathcal{K}_{\alpha}^*(\varsigma-t)}{(\varsigma-t)^{1-\alpha}}
\left(p_\omega y(x,\varsigma,u_r)-\chi_\omega^* y_d(x,\varsigma)\right)}d\varsigma\\{\kern 28pt}
=-\frac{r_2}{r_3}\frac{\mathcal{K}_{\alpha}^*(T-t)}{(T-t)^{1-\alpha}}\left(p_\omega \mathcal{M}_{\alpha}(T)y_0(x)-\chi_\omega^* y_{dT}(x)\right)\\{\kern 39 pt}
-\frac{r_1}{r_3}\int_t^T{ \frac{\mathcal{K}_{\alpha}^*(\varsigma-t)}{(\varsigma-t)^{1-\alpha}} \left(p_\omega \mathcal{M}_{\alpha}(t)y_0(x)-\chi_\omega^* y_d(x,
\varsigma)\right)}d\varsigma\\{\kern 39pt}
-\frac{r_2}{r_3}\frac{\mathcal{K}_{\alpha}^*(T-t)}{(T-t)^{1-\alpha}} p_\omega
\int_0^T{ \frac{\mathcal{K}_{\alpha}(T-\tau)}{(T-\tau)^{1-\alpha}}u_r(x,\tau)
}d\tau\\{\kern 39 pt}
-\frac{r_1}{r_3}\int_t^T{ \frac{\mathcal{K}_{\alpha}^*(\varsigma-t)}{(\varsigma-t)^{1-\alpha}}p_\omega
\int_0^\varsigma{ \frac{\mathcal{K}_{\alpha}(\varsigma-\tau)}{(\varsigma-\tau)^{1-\alpha}}
u_r(x,\tau)}d\tau}d\varsigma.
\end{array}
\end{eqnarray*}
Since  $u_r(\cdot,t)\in L^2(\Omega)$ for any given $t\in [0,T],$ let $u_{r,k}(t)=\left(u_r(\cdot,t),\xi_k\right)_{L^2(\Omega)}$. We have
\begin{eqnarray}
\begin{array}{l}
u_{r,k}(t)=\rho_k^1(t)+\rho_k^2\left(t,u_{r,k}(t)\right),
\end{array}
\end{eqnarray}
where
\begin{eqnarray*}
\begin{array}{l}
\rho_k^1(t)=-\frac{r_2}{r_3}\left(\frac{\mathcal{K}_{\alpha}^*(T-t)}{(T-t)^{1-\alpha}}\left(p_\omega \mathcal{M}_{\alpha}(T)y_0-\chi_\omega^* y_{dT}\right),
\xi_k\right)_{L^2(\Omega)}\\
-\frac{r_1}{r_3}\left(\int_t^T{ \frac{\mathcal{K}_{\alpha}^*(\varsigma-t)}{(\varsigma-t)^{1-\alpha}}
\left(p_\omega \mathcal{M}_{\alpha}(t)y_0-\chi_\omega^*  y_d(\cdot,\varsigma)\right)}d\varsigma, \xi_k\right)_{L^2(\Omega)}
\end{array}
\end{eqnarray*}
and
\begin{eqnarray*}
\begin{array}{l}
\rho_k^2\left(t,u_{r,k}(t)\right)\\
= -\frac{r_2}{r_3}\left(\frac{\mathcal{K}_{\alpha}^*(T-t)}{(T-t)^{1-\alpha}}
p_\omega\int_0^T{ \frac{\mathcal{K}_{\alpha}(T-\tau)}{(T-\tau)^{1-\alpha}}u_r(\cdot,\tau)
}d\tau,  \xi_k\right)_{L^2(\Omega)}\\
-\frac{r_1}{r_3}\left(\int_t^T{ \frac{\mathcal{K}_{\alpha}^*(\varsigma-t)}{(\varsigma-t)^{1-\alpha}}p_\omega
\int_0^\varsigma{ \frac{\mathcal{K}_{\alpha}(\varsigma-\tau)}{(\varsigma-\tau)^{1-\alpha}}
u_r(\cdot,\tau)}d\tau}d\varsigma,  \xi_k\right)_{L^2(\Omega)} .
\end{array}
\end{eqnarray*}
Then, the iterative methods can be used to obtain the values of $u_{r,k}(t)$, $t\in [0,T]$. With this, we therefore, get that the unique optimal
 solution to system $(\ref{problem})$ can be governed by
\begin{eqnarray}
\begin{array}{l}
u_r(x,t)=\sum\limits_{k\in \mathbf{N}^n}u_{r,k}(t) \xi_k(x).
\end{array}
\end{eqnarray}
This allows us to give the explicit expression of the designed optimal controller and at the same time, to minimize  the tracking
cost functional $(\ref{trackingcostf})$.

\section{Numerical example}

This section aims to present a numerical simulation illustrating our obtained results.
For the sake of simplicity, we let
$\Omega=(0,1)\subseteq \mathbf{R}$ and claim that the higher-dimensional spatial domain case
can be considered in a similar way.
\begin{figure*}
\centering
\begin{minipage}[t]{0.32\linewidth}
\includegraphics[width=1\textwidth]{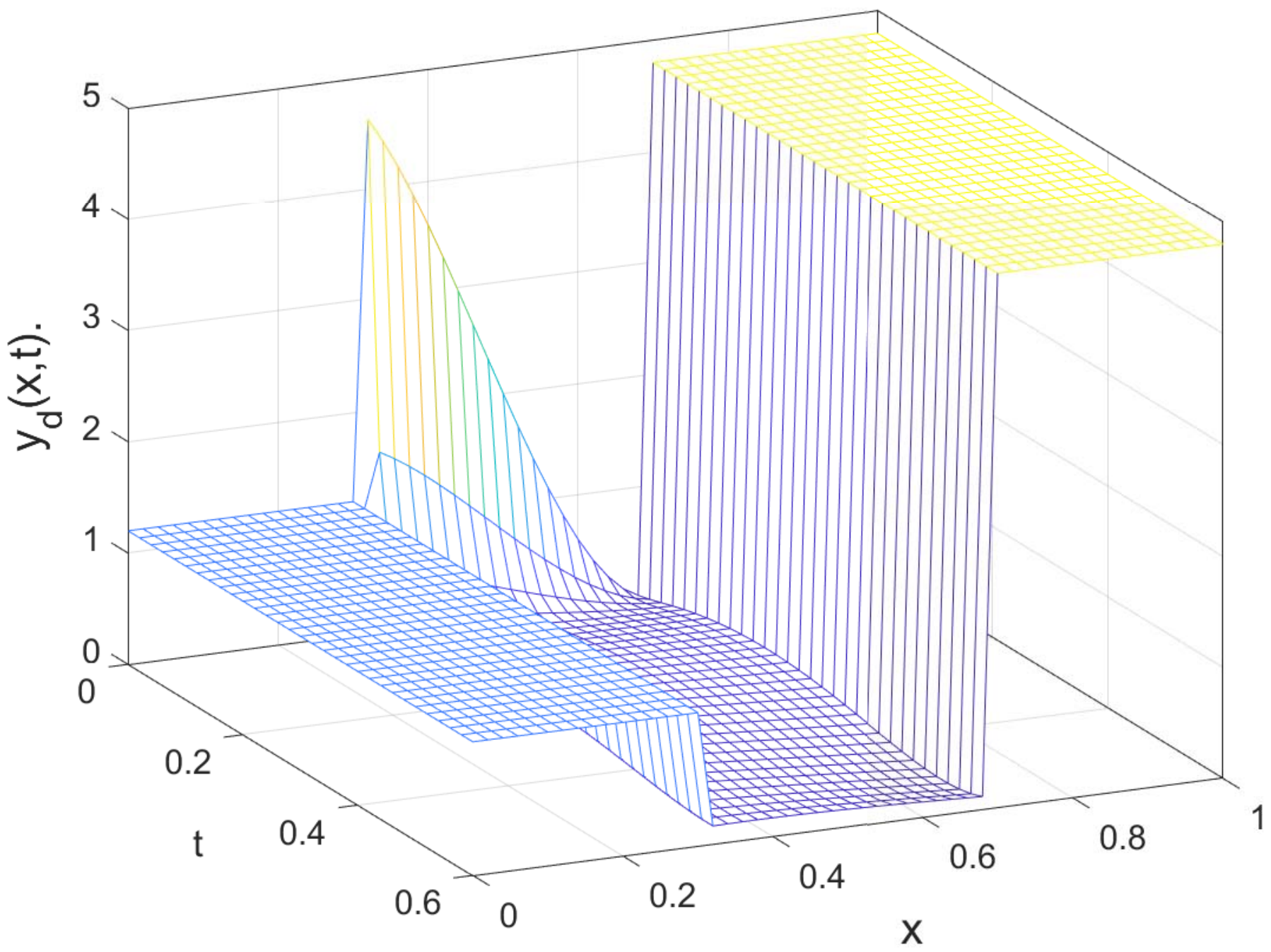}
\par{(a) The desired trajectory $y_d(x,t)$.}
\end{minipage}
\begin{minipage}[t]{0.32\linewidth}
\includegraphics[width=1\textwidth]{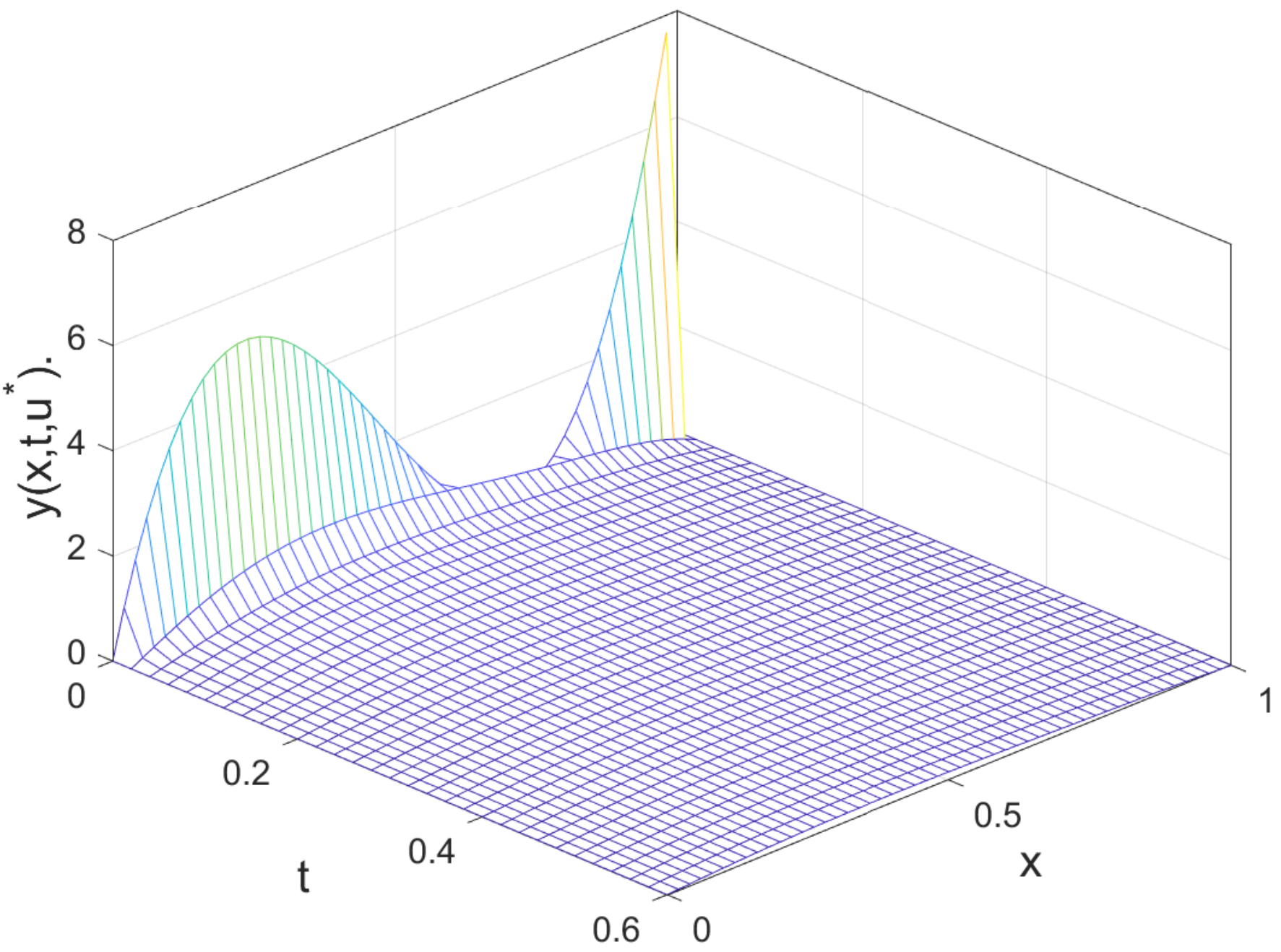}
\par{(b) The evolution of $y(x,t)$. }
\end{minipage}
\begin{minipage}[t]{0.32\linewidth}
\includegraphics[width=1\textwidth]{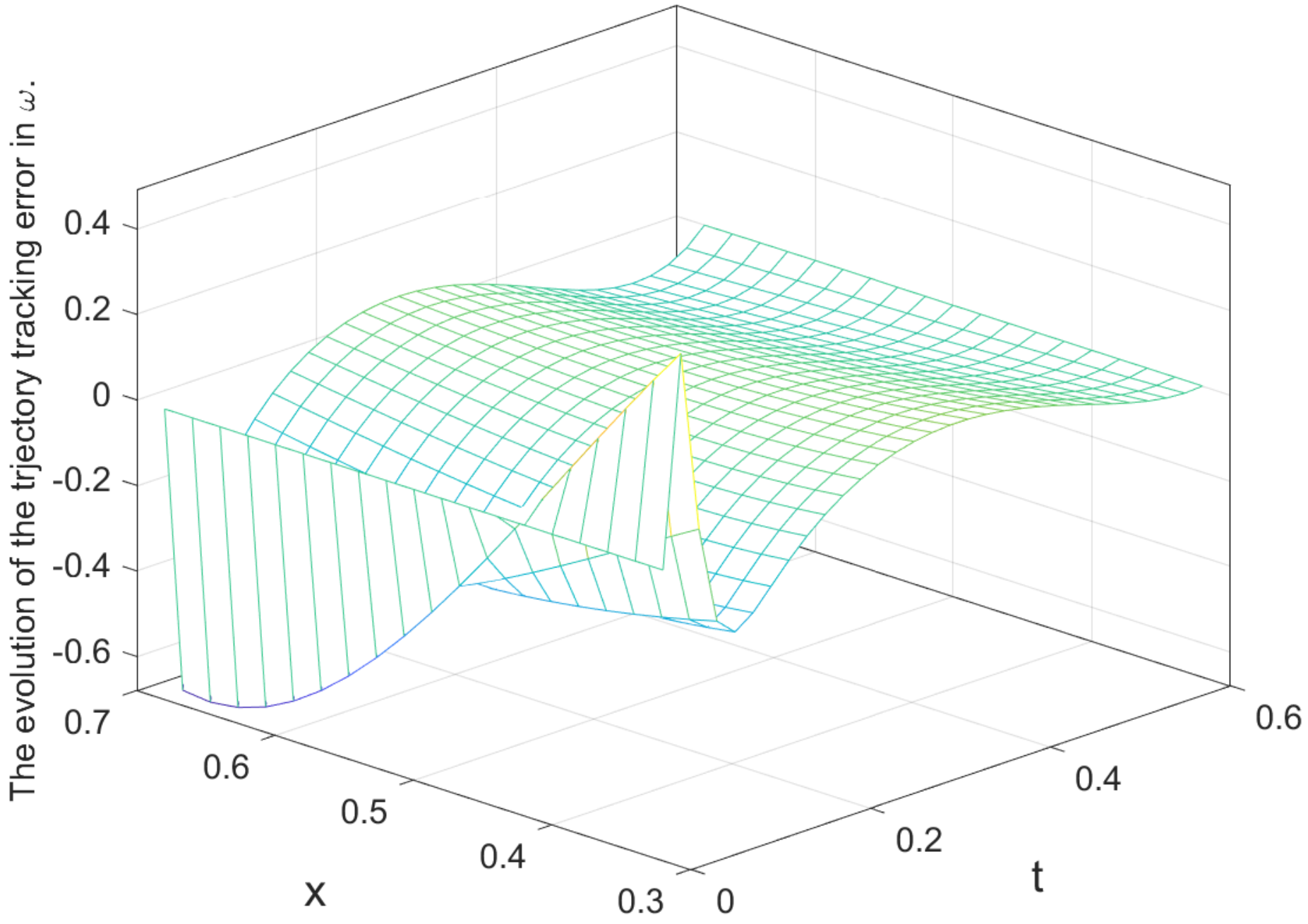}
\par{(c) The tracking error in $[0.3,0.7]$.}
\end{minipage}
\begin{minipage}[t]{0.32\linewidth}
\includegraphics[width=1\textwidth]{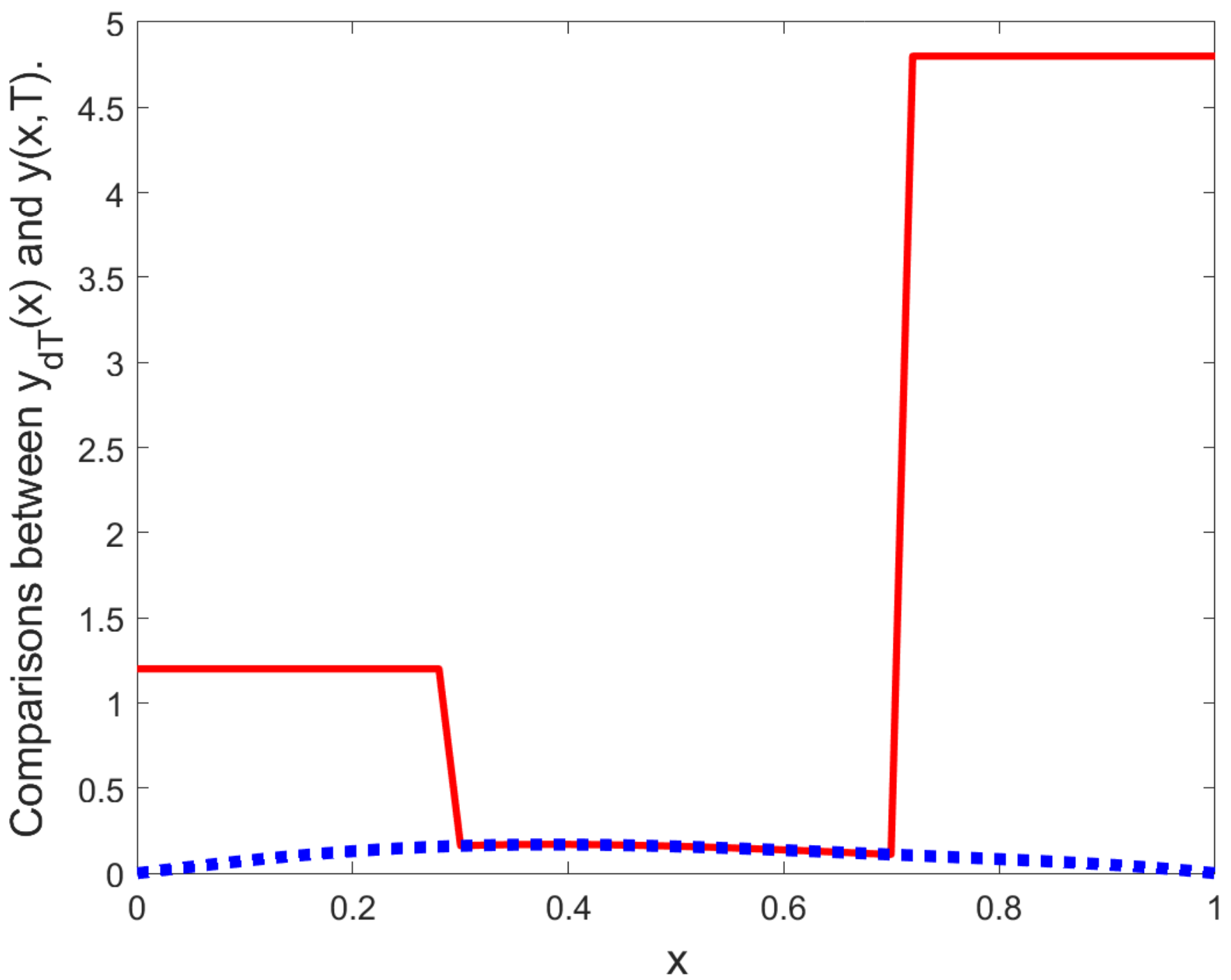}
\par{(d) Comparisons between $y_{dT}(x)$ and $y(x,T)$ in $[0.3,0.7]$.}
\end{minipage}
\hspace{0.5 cm}
\begin{minipage}[t]{0.4\linewidth}
\includegraphics[width=1\textwidth]{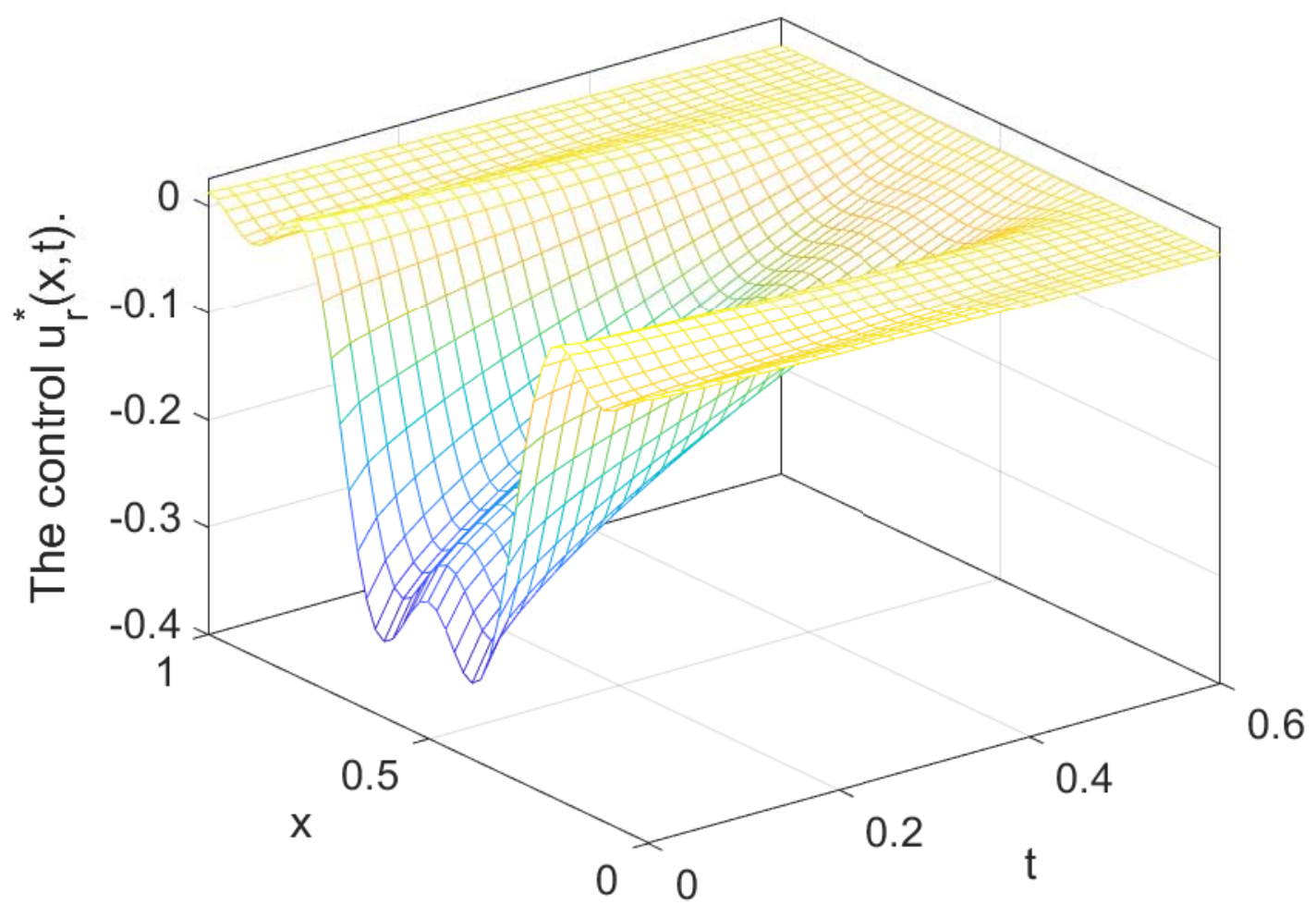}
\par{(e) The optimal control $u_r^*(x,t)$.}
\end{minipage}
\caption{Optimal regional tracking control of system $(\ref{example})$ in $[0.3,0.7]$ at time $T=0.6.$}
\label{fig2}
\end{figure*}

Let us consider  the following example
\begin{equation}\label{example}
\left\{{\begin{array}{*{20}{l}}
{}_0^CD^{0.5}_{t}y(x,t)=1.5\frac{\partial^2 y(x,t)}{\partial x^2}-y(x,t)+u(x,t)\\{\kern 51 pt}
\mbox{ in } (0,1) \times [0,0.6] ,\\
{\frac{\partial  y(0, t)}{\partial x}=\frac{\partial  y(1, t)}{\partial x}=0 \mbox{ in } [0,0.6],}\\
{ y(x,0)= 100x(x-0.7)^2\mbox{ in } (0,1).}
\end{array}}\right.
\end{equation}
Obviously, $\alpha=0.5$, $T=0.6$ and  $A=1.5\frac{\partial^2 }{\partial x^2}-1$ is a uniformly elliptic operator.
Under the Neumann boundary conditions $\frac{\partial  y(0, t)}{\partial x}=\frac{\partial  y(1, t)}{\partial x}=0$ for all
$t\in [0,0.6],$ the eigenvalue paring of operator $A$ satisfies \cite{meurer2011flatness}
\begin{eqnarray}
\begin{array}{l}
\lambda_0=-1,~\lambda_k=-1.5k^2\pi^2-1
\end{array}
\end{eqnarray}
and
\begin{eqnarray}
\begin{array}{l}
 \xi_k(x)=\left\{\begin{array}{l}
1,~~~~~~~~~~~~~~{\kern 3 pt}\mbox{if } k=0,
\\
\sqrt{2}\cos (k\pi x),~\mbox{ if } k\in \mathbf{N}\backslash \{0\}.
\end{array}\right.
\end{array}
\end{eqnarray}
Then,  the corresponding strongly continuous semigroup $\{\Phi(t)\}_{t\geq 0}$ satisfies
\begin{eqnarray}
\begin{array}{l}
\Phi(t) \varphi=\sum\limits_{k=0}^\infty{e^{\lambda_kt}\left(\varphi,\xi_k\right)_{L^2(0,1)} \xi_k},~\varphi\in L^2(0,1).
\end{array}
\end{eqnarray}
By Lemma $\ref{solution}$, it yields that
\begin{eqnarray}
\begin{array}{l}
\mathcal{M}_{\alpha}(t)\varphi= \sum\limits_{k=0}^\infty E_{\alpha}(\lambda_k t^\alpha)(\varphi,\xi_k)_{L^2(0,1)}\xi_k,~\varphi\in L^2(0,1)
\end{array}
\end{eqnarray}
and
\begin{eqnarray*}
\begin{array}{l}
\mathcal{K}_{\alpha }(t)\varphi=\sum\limits_{k=0}^\infty E_{\alpha,\alpha}
(\lambda_kt^\alpha)\left(\varphi,\xi_k\right)_{L^2(0,1)}\xi_k, ~\varphi\in L^2(0,1).
\end{array}
\end{eqnarray*}

Let the subregion  $\omega=[0.3,0.7] \subseteq (0,1)$ and the desired trajectory $y_d(x,t),$ $x\in \omega$ be
\begin{eqnarray}
\begin{array}{l}
y_d(x,t)=100x(x-0.7)^2\frac{0.6-t}{0.6e^{50t}}+4.5t(0.6-t)\\
+\frac{5t}{3}(-0.5x^4+ 2x^3-2.8\frac{x^2}{e^{0.6-t}}+1.38\frac{x}{e^{3.9-6.5t}}-0.05).
\end{array}
\end{eqnarray}
One has
\begin{eqnarray*}
\begin{array}{l}
y_{dT}(x)=y(x,0.6)\\{\kern 26pt}
=-0.5x^4+2x^3-2.8x^2+1.38x-0.05,~x\in \omega.
\end{array}
\end{eqnarray*}
In what follows, we aim to solve the following optimal control problem
\begin{eqnarray}\label{exampletrackingcostf}
\begin{array}{l}
\min\limits_{u\in L^2\left((0,1)\times [0,0.6]\right)}J(u)
\end{array}
\end{eqnarray}
with
\begin{eqnarray}
\begin{array}{l}
J(u)=10^4 \int_0^{0.6}\int_{0.3}^{0.7} {\left|\chi_{[0.3,0.7]} y(x,t,u)- y_d(x,t)\right|^2}dxdt\\{\kern 26 pt}
+10^7\int_{0.3}^{0.7}  {\left|\chi_{[0.3,0.7]}  y(x,T,u)-y_{dT}(x)\right|^2}dx
\\{\kern 26 pt}
+ \frac{1}{2}\int_{Q}{ \left|u(x,t)\right|^2} dxdt.
\end{array}
\end{eqnarray}
According to Theorem $\ref{theorem4.1}$, the optimal control problem
$(\ref{exampletrackingcostf})$ admits a unique optimal solution $u_r^*$ governed by
\begin{eqnarray}\label{u_r}
\begin{array}{l}
u_r^*(x,t)=- z(x,t,u_r^*) \mbox{ for all } (x,t)\in (0,1) \times [0,0.6].
\end{array}
\end{eqnarray}
Moreover, to illustrate the effectiveness of our results, we set
\begin{eqnarray}
h(x,t)=
\left\{\begin{array}{l}
0.25y_d(0.3,0) ~~~\mbox{ if } x\in (0,0.3),\\
y_d(x,t) ~~~~~~~~~~~\mbox{ if } x\in \omega=[0.3,0.7],\\
y_d(0.3,0)~~~~~~~~\mbox{ if } x\in (0.7,1),
\end{array}\right.\\
\begin{array}{l}
h_{dT}(x)
=\left\{\begin{array}{l}
0.25y_d(0.3,0)~~\mbox{ if } x\in (0,0.3)\\
y_{dT}(x)~~~~~~~~~~{\kern 1pt}~\mbox{if } x\in \omega=[0.3,0.7],\\
y_d(0.3,0)~~~~~~~\mbox{ if } x\in (0.7,1)
\end{array}
\right.
\end{array}
\end{eqnarray}
and then plot Figures $\ref{fig2}$ $(a)-(d)$, which show how closely does the  state evolution track along the desired
trajectory $y_d(x,t)\in L^2\left((0.3,0.7) \times [0,0.6]\right)$ and the final state reach the  target function $y_{dT}=y_d(\cdot,T)
\in L^2(0.3,0.7)$ in $\omega=[0.3,0.7]$ with the error
\begin{eqnarray}
\begin{array}{l}
\left\|\chi_{(0.3,0.7)} y(\cdot,T) -y_T (x)\right\|_{L^2(0.3,0.7)}\leqslant 0.005,\\
{\kern 1 pt}\left\|\chi_{(0.3,0.7)}y(\cdot,t) -y_d(x,t) \right\|_{L^2(0.3,0.7)}\leqslant 2.99.
\end{array}
\end{eqnarray}
The corresponding optimal solution of the optimal control problem $(\ref{exampletrackingcostf})$ is depicted
in Figure $\ref{fig2}$ $(e)$ with the costs $\|u_r^*\|_{L^2(0,1)\times L^2(0,0.6)}=12.57$.

\section{Conclusion}
Sufficient and necessary conditions for optimal regional trajectory tracking control problem of linear time-fractional
diffusion systems are obtained in this paper by using the HUM. The obtained results not only can be used directly to
discuss the systems that are not controllable on the whole domain, but also yield an explicit expression of the control
signal in terms of the desired trajectory and minimize the proposed tracking cost functional as well. This is very appealing
in practical applications and pose many new theoretically challenges at the same time. Moreover, we claim that the main results
in this paper can be extended to more complex fractional-order distributed parameter systems (see those in \cite{ge2015JAS} for
example)and various open questions such as optimal actuation configuration problems for regional tracking control of the coupled
nonlinear space-time fractional diffusion systems are still open.

\bibliographystyle{elsarticle-num}
\bibliography{PDEs}

\end{document}